\documentstyle[11pt]{article}
\input amssym.def
\input amssym

\def\P{\Bbb P}

\newcommand{\pf}{\noindent{\it Proof.}\ }
\newcommand{\qed}{\hfill $\Box$ }
\def\sect#1{\section*{{\large\bf  #1}}}
\def\jielun#1{{\par\noindent\hskip 0.5 cm}
{\bf #1}\bgroup\it}
\def\endjielun{\egroup\par\bigskip}
\def\bdefi#1{{\par\noindent\hskip 0.5 cm}
{\bf #1}\bgroup\rm}
\def\edefi{\egroup\par\bigskip}
\def\bli#1{{\par\noindent\hskip 0.5 cm}
{\bf #1}\bgroup\rm}
\def\eli{\egroup\par\bigskip}
\begin{document}
\title{\bf   {\CC {\char64}n K+ 4z
J{\char125} Sk A? WS 4z J{\char125} 5D JU Ku}}
\title{{\bf {A Note on
Anti-Pluricanonical Maps for 5-Folds}}
\thanks{Research partially supported by SFB-237 of the DFG.}
\thanks{Keywords: anti-pluricanonical map,  birationality.}
\thanks{\it Mathematics Subject Classifications(2000):
 \ {\rm 14E05, 14J30, 14N05.}}}
\author{Qi-Lin Yang\\
Department of Mathematics,
Tsinghua University \\
Beijing 100084, P. R. China\\
(E-mail: qlyang@math.tsinghua.edu.cn)\\
  }
\maketitle
\date{}

\begin{abstract}
We prove that the anti-pluricanonical map $\Phi_{|-mK_{X}|}$ is
birational when $m\geq 16$ for $5$-fold $X$ whose anticanonical
divisor is nef and big.
\end{abstract}

\sect{\S1. \  Introduction.}
 Throughout the ground field $k$ is
always supposed to be algebraically closed of characteristic zero.
 Let $X$ be a non-singular  $n$-fold over $k$ and assume
its anticanonical
 divisor $-K_{X}$ is a nef and big divisor. It is an interesting problem
 to find an explicit lower bound $l(n)$ such that the rational map
 $\Phi _{|-mK_{X}|}$ associated with the
 complete linear system $|-mK_{X}|$ is a birational map
 onto its image for any $m\geq l(n).$ Ando ([1, Theorem 9]) first
 gave the bounds $l(2)=3, l(3)=5$ and $l(4)=12.$ Fukuda [3] improved
  Ando's method and got the bounds $l(3)=4,  l(4)=11$ and
  $$l(n)=2^{n-2}\cdot (n+4[n/2]-1)-2[n/2]-1$$
  for any $n\geq 5.$ Chen [2] also used the similar ideas to improve
   Ando's partial results.
   Using  the  Key Lemma  in [3], this note proves

\bli{Main Theorem} Let $X$ be a smooth 5-fold whose anticanonical
divisor $-K_{X}$ is nef and big. Then $\Phi _{|-mK_{X}|}$ is a
birational map when $m\geq 16.$ \eli

\sect{\S2. \  Preparations.}
 In this note we use the standard terminology as in [4, 5].
 For example, $c_{i}:=c_{i}(TX)$ is the $i$-th Chern
 class of the tangential bundle; $H^{i}(X, {\cal F})$
  denotes the $i$-th cohomolgy with coefficient in a
  coherent sheaf ${\cal F},$ and
  $h^{i}(X, {\cal F})=\dim _{k}H^{i}(X, {\cal F}).$
We simply denote $H^i (X, {\cal O}_{X}(D))$  by $H^i(X,D)$ if the
sheaf is induced by a divisor $D.$

  We will use Lemma 1, a special case of the Key Lemma in [3], which improved the
  Theorem 5  in  [1].
\bli{Lemma 1} (Ando [1], Fukuda [3],  Chen [2]). Let $X$ be a
nonsingular projective variety of dimension $n$ and  $-K_{X}$ is a
nef and big divisor. We assume:

(i)   For each $i$ with $1\leq i\leq n-2,$ there exists a natural
number $r_i$ such that  $\dim \Phi_{|-r_{i}K_{X}|} (X)\geq i.$

(ii) There exist an integer $r_{0}\geq 3$ such that $H^{0}(X,
-rK_{X})\not=0$ for any $r\geq r_{0}.$

Then $\Phi_{|-mK_{X}|}$ is birational for all $m\geq
r_{0}+(r_{1}+\cdots r_{n-2}).$ \eli

\pf In the Key Lemma in [3],  we let the nef and big divisor
$R=-K_{X}$ and the numerically trivial divisor $T=0.$ By our
assumptions we have $H^0 (X, -rK_{X})=H^0 (X,
-(r+1)K_{X}+K_{X})\not=0$ for any $r+1\geq
r_{0}+1:=\hat{r}_{0}\geq 4.$ So both (1) and (2) of the Key Lemma
are satisfied. Hence $\Phi_{|-mK_{X}|}=\Phi_{|-(m+1)K_{X}+K_{X}|}$
is birational when $m+1\geq \hat{r}_{0}+(r_{1}+\cdots r_{n-2}),$
thus $\Phi_{|-mK_{X}|}$ is birational for all $m\geq
r_{0}+(r_{1}+\cdots r_{n-2}).$\qed

To use Lemma 1, we need the following Lemma, it is the Proposition
6 in [1], we refer to [1, 2, 3] for reference of it.

 \bli{Lemma 2} (Matsusaka \& Maehara). Let $D$ be a nef and big
divisor and $\dim X=n.$ If $h^{0}(X, mD)>m^{r}D^n +r,$ then  $\dim
\Phi_{|mD|}>r.$
 \eli

\sect{\S3. \  Proof of the Main Theorem.}
 Let
    $$P(m):=\chi ({\cal O}_{X}(-mK_{X}) )=
    \sum_{i=0}^{5}h^{i}(X, -mK_{X})=
h^{0}(X,-mK_{X}),$$
    since $-K_{X}$ is nef and big, by
    the Kawamata-Viehweg vanishing theorem
     ({\rm cf.} [5, Corollary 1-2-2])
     we have
      $H^{i}(X, -mK_{X})=0$ for $m\geq 0$ and $i>0.$
    Thus  $\chi ({\cal O}_{X})=1.$
    Note by definition,  $c_{1}=-K_{X},$
    and $c_{i}=0$ for $i>5.$ Combine these facts with Hirzebruch-Riemann-Roch
    formula ([4, $P_{432}$]), we have
      $$\begin{array}{rcl}
      P(m)&=&\int_{X}{\rm ch}({\cal O}_{X}(-mK_{X})){\rm Td}(X)\\
         &=& m(m+1)(2m+1)(3m^2 + 3m -1)\frac{(-K_{X})^{5}}{720}
          +m(m+1)(2m+1)\frac{(-K_{X})^{3}\cdot c_{2}}{144}+(2m+1)\\
         &=&(2m+1)\{(m(m+1)[(3m^2 + 3m -1)a+b]+1\},
      \end{array}$$
     where $720a=(-K_{X})^{5}$ and
     $144b=(-K_{X})^{3}\cdot c_{2}.$

To use Lemma 1 and Lemma 2, we need  priori estimates of $P(m)$
for $m\geq 0.$
 \bli{Proposition 1}.
$$\begin{array}{rcl}
{\rm (i)}&& P(1)=0, P(2)\geq 0 \Rightarrow P(3)\geq 35;\\
{\rm (ii)} &&P(1)=1, P(2)\geq 1 \Rightarrow P(3)\geq 21;\\
{\rm (iii)}&&P(1)=2,P(2)\geq 2 \Rightarrow P(3)\geq 7;\\
 {\rm (iv)}&& P(1)=3\Rightarrow P(2)\geq 6;\\
{\rm (v)} &&P(1)=3, P(2)=6 \Rightarrow P(3)= 49;\\
{\rm (vi)} &&P(m+1)>P(m) ~~{\rm when}~~m>3~~{\rm and}~~P(3)\geq 7.\\
 \end{array}$$
 \eli
\pf  Assume $P(1)=3[2(5a+b)+1]:=l\geq 0,$ we have
$b=\frac{1}{6}(l-3)-5a.$ By $P(2)\geq P(1)$ we get $a\geq
\frac{1}{360}(10-4l),$  so we have $P(3)\geq 35-14l.$ Hence we get
(i)-(iii).

Now we assume $P(1)=3$ and $P(2)=5[6(17a+b)+1]:=l.$ Then $b=-5a$
and $l=5(12a+1)>5,$ and hence $P(2)\geq 6$ and we have (iv). If
$P(2)=2P(1)=6.$ Then we have $a=\frac{1}{60}$ and
$b=-\frac{1}{12},$ So $P(3)=7[12(35a+b)+1]=49>7$ we get (v).

If $P(2)>6,$ then $P(3)\geq P(2)\geq 7.$ Combine with (i)-(v) we
always have $P(3)\geq 7.$ Since $P(1)\geq 0,$ we have $b\geq
-5a-\frac{1}{2}$ and $[3m(m+1)-1]a+b>0$ when $m\geq 3.$ Thus
$P(m+1)-P(m)>(m+1)\{[3m(m+1)-1]a+b+1\}[(m+2)(2m+3)-(m+1)(2m+1)]>0$
when $m\geq 3,$ we get (vi). \qed
 \bli{Proposition 2}.
$$\begin{array}{rcl}
{\rm (i)} &&\dim \Phi_{|-mK_{X}|}(X)\geq 1, ~~{\rm for~~ any}~~ m\geq 3;\\
 {\rm (ii)} &&\dim \Phi_{|-mK_{X}|}(X)\geq 2, ~~{\rm for~~ any}~~ m\geq 4;\\
{\rm (iii)}&&\dim \Phi_{|-mK_{X}|}(X)\geq 3,  ~~{\rm for~~any}~~ m\geq 6.\\
  \end{array}$$
 \eli

\pf By Proposition 1, $P(3)=7[12(35a+b)+1]\geq 7,$ so we have (i)
and $b\geq -35a.$ Thus $P(4)=9[20(59a+b)+1]\geq 180\times 24 a+9
>6(-K_{X})^{5}+2.$ By Lemma 2 we have (ii).
$P(6)=13[42(125a+b)+1]\geq \frac{13\times
21}{4}(-K_{X})^{5}+13>36(-K_{X})^{5}+3,$ we have (iii). \qed

\vskip 0.5cm

{\it Proof of Main Theorem.} By Proposition 1 we have  $h^{0}(X,
-3K_{X})\geq 7,$ so we can put $r_{0}=3.$ By Proposition 2, we can
set $r_{1}=3,r_{2}=4,r_{3}=6.$ By Lemma 1 when $m\geq
r_{0}+r_{2}+r_{3}+r_{4}=16,$ then $\Phi_{|-mK_{X}|}$ is a
birational map.\qed

 \sect{\S4. \ An example.}
\bli{Example 1}. Let $\pi: E={\cal O}_{{\P}^{1}} \oplus {\cal
O}_{{\P}^{1}}\oplus {\cal O}_{{\P}^{1}}\oplus {\cal
O}_{{\P}^{1}}\oplus {\cal O}_{{\P}^{1}}(1)\rightarrow {\P}^{1}$ be
a rank $5$ vector bundle. Let $X={I\!\!P}(E).$ Then by
calculations in the Exercise 8.4 of [4, $P_{253}$],
$K_{X}=-5L+\pi^{*}(\det (E)+K_{{\P}^{1}})=-5L-H,$ where $L\in
|{\cal O}_{X}(1)|, H\in |\pi^{*}{\cal O}_{{\P}^{1}}(1)|.$ Clearly
$X$ is a Fano manifold and $-K_X$ is a nef and big divisor. Note
$L^{5}=H\cdot L^{4}=1,$ so $(-K_{X})^{5}=2\times 5^{5}.$ By Leray
spectral sequence and the fact that $R^{i}\pi_{*}({\cal O}(l))=0$
for any $i>0$ and $l>-5,$ we have
 $H^{i}(X,-mK_{X})=0$ when $i>0$ and
 $H^{0}(X, -mK_{X})=H^{0}(X,{\cal O}_{X}(5m)\otimes \pi^{*}
 {\cal O}_{{\P}^{1}}(m))= H^{0}({\P}^{1}, S^{5m}(E) \otimes
  {\cal O}_{{\P}^{1}}(m))$ for any $m>0.$  Note that
$$\begin{array}{rcl} S^{5m}(E)&=&\bigoplus_{i=0}^{5m}
S^{5m-i}({\cal O}_{{\P}^{1}}\oplus {\cal O}_{{\P}^{1}}\oplus {\cal
O}_{{\P}^{1}}\oplus {\cal O}_{{\P}^{1}})
\otimes {\cal O}_{{\P}^{1}}(i) \\
&=&\bigoplus_{i=0}^{5m}({{\cal O}_{{\P}^{1}}(i)\oplus {\cal
O}_{{\P}^{1}}(i)\oplus
\cdots \oplus {\cal O}_{{\P}^{1}}(i)}),\\
\end{array}$$
it is a bundle of rank $\frac{1}{6}(5m-i-1)(5m-i)(5m-i+1)$ in the
last bracket of above summation. So,
$$\begin{array}{rcl}
h^{0} ({\P}^{1}, S^{5m}(E) \otimes {\cal
O}_{{\P}^{1}}(m))&=&\frac{1}{6}\sum_{i=0}^{5m}(5m-i-1)(5m-i)(5m-i+1)
h^{0}({\P}^{1}, {\cal O}_{{\P}^{1}}(m+i))\\
 &=&\frac{1}{6}\sum_{i=0}^{5m}(m+i+1)(5m-i-1)(5m-i)(5m-i+1).\\
 &=&\frac{1}{24}m(5m-1)(5m+1)(5m+2)(10m+3).
\end{array}$$
It is easy to check that $h^{0}({I\!\!P}(E), -K_{X})=91,$ we can
take $r_{0}=r_1 =3;$ and
 $h^{0}({I\!\!P}(E),-4K_{X})=62909> 10(-K_{X})^{5}+2,$ we take $r_{2}=4,$
 and $h^{0}({I\!\!P}(E),-5K_{X})=186030> 5^{2}(-K_{X})^{5}+3$ we can
  take $r_{3}=5.$  So
  $\Phi _{|-mK_{X}|}$ is a birational map when $m\geq 15.$
\eli

 {\bf Question}. Find out the lowest bound $l(n)$ such that
$\Phi_{|-mK_{X}|}$ is  birational when $m\geq l(n).$

 We also don't know  how to improve the bounds  $l(n)$ given in [3]
 for $n>5$ since the Hirzebruch-Riemann-Roch
    formula is more complicate in these cases.
\section*{Acknowledgement}

The author would like to thank the referee for showing him the
paper [3], base on which he improve the earlier bound to present
stage.
 Part of this work was done while I was
 visiting as a guest fellow at
 the Institut f\"ur Mathematik, Ruhr Universit\"at Bochum, Germany.
I would like to thank Prof. A. Huckleberry and P. Heinzner
 for showing me the thesis
of Dr. S. Kebekus, which stimulates my interest in algebraic
geometry.

\begin{center}
{\bf\Large References}
\end{center}
{\parindent=0pt
\def\toto#1#2{\centerline{\makebox[1.5cm][l] {#1\hss}
\parbox[t]{13cm}{#2}}\vspace{\baselineskip}}

\toto{[1]}{{\rm T. ANDO,} { Pluricanonical systems of algebraic
varieties of general type of dimension $\leq$ 5,} Adv. Stud. in
Pure Math., {\bf 10}(1987), {\it Algebraic Geometry, Sendi, 1985,}
1-10.}

\toto{[2]}{{\rm M. CHEN,} {A note on pluricanonical maps for
varieties of dimension 4 and 5,}  J. Math. Kyoto Univ., {\bf
37}(1997) 513-517.}

\toto{[3]}{{\rm S. FUKUDA,} {A note on Ando's paper
``pluricanonical systems of algebraic varieties of general type of
dimension $\leq$ 5,}  Tokyo J. Math., {\bf 14}(1991) 479-487.}

\toto{[4]}{{\rm R. HARTSHORNE,} {\it Algebraic Geometry,}
  GTM {\bf 52}, Springer-Verlag, 1977.}

\toto{[5]}{{\rm  Y. KAWAMATA, K.  MATSUDA, \& K. MATSUKI,} {
Introduction to minimal model problem,} Adv. Stud. in Pure Math.,
{\bf 10}(1987), {\it Algebraic Geometry, Sendi, 1985,} 283-360.}

\end{document}